\input amstex
\input amsppt.sty
\magnification=\magstep1 \hsize=32truecc
 \vsize=22.5truecm
\baselineskip=14truept
\NoBlackBoxes
\def\q{\quad}
\def\qq{\qquad}
\def\mod#1{\ (\text{\rm mod}\ #1)}
\def\t{\text}
\def\qtq#1{\q\t{#1}\q}
\def\f{\frac}
\def\e{\equiv}
\def\a{\alpha}
\def\b{\binom}

\def\sls#1#2{(\f{#1}{#2})}
 \def\ls#1#2{\big(\f{#1}{#2}\big)}
\def\Ls#1#2{\Big(\f{#1}{#2}\Big)}

\let \pro=\proclaim
\let \endpro=\endproclaim

\topmatter
\title\nofrills Some conjectures on congruences
\endtitle
\author ZHI-Hong Sun\endauthor
\affil School of Mathematical Sciences, Huaiyin Normal University,
\\ Huaian, Jiangsu 223001, PR China
\\ Email: zhihongsun$\@$yahoo.com
\\ Homepage: http://www.hytc.edu.cn/xsjl/szh
\endaffil

 \nologo \NoRunningHeads

\abstract{Let $p$ be an odd prime. In the paper we collect the
author's various conjectures on congruences modulo $p$ or $p^2$,
which are concerned with sums of binomial coefficients, Lucas
sequences, power residues and special binary quadratic forms.
\par\q
 }
 \endabstract

\endtopmatter
\document
\subheading{1. Notation}
\par Let $\Bbb Z$ and $\Bbb N$ be the sets
of integers and positive integers respectively. For $b,c\in\Bbb Z$
the Lucas sequences $\{U_n(b,c)\}$ and $\{V_n(b,c)\}$ are defined by
$$\aligned &U_0(b,c)=0,\ U_1(b,c)=1,\\& U_{n+1}(b,c)=bU_n(b,c)-cU_{n-1}(b,c)\
(n\ge 1)\endaligned\tag 1.1$$ and $$\aligned &V_0(b,c)=2,\
V_1(b,c)=b,\\& V_{n+1}(b,c)=bV_n(b,c)-cV_{n-1}(b,c)\ (n\ge
1).\endaligned\tag 1.2$$ Let $d=b^2-4c$. It is well known that for
$n\in\Bbb N$,
$$U_n(b,c)=\cases \f 1{\sqrt d}\big\{\ls{b+\sqrt d}2^n-
\ls{b-\sqrt d}2^n\big\}&\t{if $d\not=0$,}
\\n(\f b2)^{n-1}&\t{if $d=0$}\endcases\tag 1.3$$
and
$$V_n(b,c)=\Ls {b+\sqrt d}2^n+
\Ls{b-\sqrt d}2^n.\tag 1.4$$
\par  Let $[x]$ be the greatest integer not exceeding $x$,
and let $\sls am$ be the Jacobi symbol. For $m\in\Bbb Z$ with
$m=2^{\alpha}m_0(2\nmid m_0)$ we say that $2^{\alpha} \ \Vert\ m$
and $m_0$ is the odd part of $m$. For $t\in\Bbb Z$ let
$$\delta(t)=\cases 1&\t{if $8\mid t$,}
\\-1&\t{if $8\nmid t$.}
\endcases$$
\par For an integer $m$ and odd prime $p$ with $p\nmid m$ let
$$Z_p(m)=\sum_{n=0}^{p-1}\f{\b{2n}n}{m^n}\sum_{k=0}^n\b nk^3.$$

 \subheading{2. Conjectures on  power residues}

\par In 1980 and 1984 Hudson and Williams proved the following
result.
 \pro{Theorem 2.1} Let $p\e
1\mod{24}$ be a prime and hence $p=c^2+d^2=x^2+3y^2$ for some
$c,d,x,y\in\Bbb Z$. Suppose $c\e 1\mod 4$.
\par $(\t{\rm i) ([HW]})$ If $c\e \pm (-1)^{\f y4}\mod 3$, then
$3^{\f{p-1}8}\e \pm 1\mod p$.
\par $(\t{\rm ii) ([H]})$ If $d\e \pm (-1)^{\f y4}\mod 3$, then
$3^{\f{p-1}8}\e \pm \f dc\mod p$.
\endpro
\par Hudson and Williams proved Theorem 2.1(i) by using the
cyclotomic  numbers of order $12$, and Hudson proved Theorem 2.1(ii)
using the Jacobi sums of order $24$.
\par In [S3] the author posed the following conjectures similar to Theorem 2.1.
\pro{Conjecture 2.1 ([S3, Conjecture 9.1])} Let $p\e 13\mod{24}$ be
a prime and hence $p=c^2+d^2=x^2+3y^2$ for some $c,d,x,y\in\Bbb Z$.
Suppose $c\e 1\mod 4$, $x=2^{\a}x_0$, $y=2^{\beta}y_0$ and $x_0\e
y_0\e 1\mod 4$. Then
$$3^{\f{p-5}8}\e\cases \pm\f yx\mod p&\t{if $x\e \pm c\mod 3$,}
\\\mp\f{dy}{cx}\mod p&\t{if $x\e\pm d\mod 3$.}
\endcases$$\endpro
\par Conjecture 2.1 has been checked for all primes $p<20,000$.
 \pro{Conjecture 2.2 ([S3, Conjecture 9.2])} Let $p\e 1,9,25\mod{28}$ be a
prime and hence $p=c^2+d^2=x^2+7y^2$ for some $c,d,x,y\in\Bbb Z$.
Suppose $c\e 1\mod 4$, $x=2^{\a}x_0$, $y=2^{\beta}y_0$ and $x_0\e
y_0\e 1\mod 4$.
\par $(\t{\rm i})$ If $p\e 1\mod 8$, then
$$7^{\f{p-1}8}\e \cases -(-1)^{\f y4}\mod p&\t{if $7\mid c$,}
\\(-1)^{\f y4}\mod p&\t{if $7\mid d$,}
\\\mp (-1)^{\f y4}\f dc\mod p&\t{if $c\e\pm d\mod 7$.}
\endcases$$
\par $(\t{\rm ii})$ If $p\e 5\mod 8$, then
$$7^{\f{p-5}8}\e \cases -\f yx\mod p&\t{if $7\mid c$,}
\\\f yx\mod p&\t{if $7\mid d$,}
\\\mp \f {dy}{cx}\mod p&\t{if $c\e\pm d\mod 7$.}
\endcases$$
\endpro
\par Conjecture 2.2 has been checked for all primes $p<20,000$.
\pro{Conjecture 2.3 ([S3, Conjecture 9.7])} Let $p\e 1,9\mod{20}$ be
a prime and hence $p=c^2+d^2=x^2+5y^2$ for some $c,d,x,y\in\Bbb Z$.
Suppose $c\e 1\mod 4$ and all the odd parts of $d,x,y$ are congruent
to $1$ modulo $4$.
\par $(\t{\rm i})$ If $p\e 1\mod 8$, then
$$5^{\f{p-1}8}\e \cases \pm (-1)^{\f d4}\delta(y)\mod p
&\t{if $x\e \pm c\mod 5$,}
\\\pm(-1)^{\f d4}\delta(y)\f dc\mod p
&\t{if $x\e \pm d\mod 5$.}\endcases$$
\par $(\t{\rm ii})$ If $p\e 5\mod 8$, then
$$5^{\f{p-5}8}\e \cases \pm \delta(x)\f{dy}{cx}
\mod p &\t{if $x\e \pm c\mod 5$,}
\\\mp \delta(x)\f yx\mod p
&\t{if $x\e \pm d\mod 5$.}\endcases$$
\endpro
\par Conjecture 2.3 has been checked for all primes $p<20,000$.

\pro{Conjecture 2.4 ([S3, Conjecture 9.8])} Let $p\e 1,9\mod {40}$
be a prime and hence $p=c^2+d^2=x^2+10y^2$ for some $c,d,x,y\in\Bbb
Z$. Suppose $c\e x\e 1\mod 4$. Then
$$5^{\f{p-1}8}\e \cases\pm (-1)^{\f d4+\f {x-1}4}\f dc\mod p&\t{if $x\e \pm
d\mod 5$,}\\\pm (-1)^{\f d4+\f {x-1}4}\mod p&\t{if $x\e \pm c\mod
5$.}
\endcases$$\endpro

 \pro{Conjecture 2.5 ([S3, Conjecture 9.9])} Let $p\e 1,9,17,25,29,49\mod {52}$ be a prime and
hence $p=c^2+d^2=x^2+13y^2$ for some $c,d,x,y\in\Bbb Z$. Suppose
$c\e 1\mod 4$ and all the odd parts of $d,x,y$ are congruent to $1$
modulo $4$.
\par $(\t{\rm i})$ If $p\e 1\mod 8$, then
$$13^{\f{p-1}8}\e \cases \mp (-1)^{\f d4}\delta(y)\f dc\mod p
&\t{if $\f{2c+3d}x\e \pm 1,\pm 3,\pm 9\mod {13}$,}
\\\pm(-1)^{\f d4}\delta(y)\mod p
&\t{if $\f{2c+3d}x\e \pm 2,\pm 5,\pm 6\mod {13}$.}\endcases$$
\par $(\t{\rm ii})$ If $p\e 5\mod 8$, then
$$13^{\f{p-5}8}\e \cases \pm \delta(x)\f yx
\mod p &\t{if $\f{2c+3d}x\e \pm 1,\pm 3,\pm 9\mod {13}$,}
\\\pm \delta(x)\f {dy}{cx}\mod p
&\t{if $\f{2c+3d}x\e \pm 2,\pm 5,\pm 6\mod {13}$.}\endcases$$
\endpro

\pro{Conjecture 2.6 ([S3, Conjecture 9.16])} Let $p\e 1\mod 4$ be a
prime and $p=c^2+d^2=x^2+17y^2$ for some $c,d,x,y\in\Bbb Z$. Suppose
$c\e 1\mod 4$ and all the odd parts of $d,x,y$ are numbers of the
form $4k+1$.
\par $(\t{\rm i})$ If $p\e 1\mod 8$, then
$$17^{\f{p-1}8}\e \cases
-(-1)^{\f d4+\f{xy}4}\f dc\mod p&\t{if $4c+d\e \pm 6x,\pm
7x\mod{17}$,}
\\(-1)^{\f d4+\f{xy}4}\f dc\mod p&\t{if $4c+d\e \pm 3x,\pm 5x\mod{17}$,}
\\(-1)^{\f d4+\f{xy}4}\mod p&\t{if $4c+d\e \pm x,\pm 4x\mod{17}$,}
\\-(-1)^{\f d4+\f{xy}4}\mod p&\t{if $4c+d\e \pm 2x,\pm 8x\mod{17}$.}
\endcases$$
\par $(\t{\rm ii})$ If $p\e 5\mod 8$, then
$$17^{\f{p-5}8}\e \cases
(-1)^x\f yx\mod p&\t{if $4c+d\e \pm 6x,\pm 7x\mod{17}$,}
\\-(-1)^x\f yx\mod p&\t{if $4c+d\e \pm 3x,\pm 5x\mod{17}$,}
\\(-1)^x\f {dy}{cx}\mod p&\t{if $4c+d\e \pm x,\pm 4x\mod{17}$,}
\\-(-1)^x\f {dy}{cx}\mod p&\t{if $4c+d\e \pm 2x,\pm 8x\mod{17}$.}
\endcases$$
\endpro

\pro{Conjecture 2.7 ([S5, Conjecture 4.3])} Let $p\e 1\mod 4$ and
$q\e 3\mod{8}$ be primes such that $p=c^2+d^2=x^2+qy^2$ with
$c,d,x,y\in\Bbb Z$ and $q\mid cd$. Suppose $c\e 1\mod 4$,
$x=2^{\a}x_0$, $y=2^{\beta}y_0$ and $x_0\e y_0\e 1\mod 4$.
\par $(\t{\rm i})$ If $p\e 1\mod 8$, then
$$q^{\f{p-1}8}\e \cases \pm (-1)^{\f y4}\mod p&\t{if $x\e \pm c\mod{q}$,}
\\\mp(-1)^{\f{q-3}8+\f y4}\f dc\mod p&\t{if $x\e\pm
d\mod{q}$.}
\endcases$$
\par $(\t{\rm ii})$ If $p\e 5\mod 8$, then
$$q^{\f{p-5}8}\e \cases \pm \f yx\mod p&\t{if $x\e\pm
c\mod{q}$,}\\\mp (-1)^{\f{q-3}8}\f{dy}{cx}\mod p&\t{if $x\e \pm
d\mod{q}$.}
\endcases$$
\endpro

\pro{Conjecture 2.8 ([S5, Conjecture 4.4])} Let $p\e 1\mod 4$ and
$q\e 7\mod{16}$ be primes such that $p=c^2+d^2=x^2+qy^2$ with
$c,d,x,y\in\Bbb Z$ and $q\mid cd$. Suppose $c\e 1\mod 4$,
$x=2^{\a}x_0$, $y=2^{\beta}y_0$ and $x_0\e y_0\e 1\mod 4$.
\par $(\t{\rm i})$ If $p\e 1\mod 8$, then
$$q^{\f{p-1}8}\e \cases (-1)^{\f y4}\mod p&\t{if $q\mid d$,}
\\-(-1)^{\f y4}\mod p&\t{if $q\mid c$.}
\endcases$$
\par $(\t{\rm ii})$ If $p\e 5\mod 8$, then
$$q^{\f{p-5}8}\e \cases \f yx\mod p&\t{if $q\mid d$,}
\\-\f yx\mod p&\t{if $q\mid c$.}
\endcases$$
\endpro

\pro{Conjecture 2.9([S5, Conjecture 4.5])} Let $p\e 1\mod 4$ and
$q\e 15\mod{16}$ be primes such that $p=c^2+d^2=x^2+qy^2$ with
$c,d,x,y\in\Bbb Z$ and $q\mid cd$. Suppose $c\e 1\mod 4$,
$x=2^{\a}x_0$, $y=2^{\beta}y_0$ and $x_0\e y_0\e 1\mod 4$.
\par $(\t{\rm i})$ If $p\e 1\mod 8$, then
$q^{\f{p-1}8}\e (-1)^{\f y4}\mod p.$
\par $(\t{\rm ii})$ If $p\e 5\mod 8$, then
$q^{\f{p-5}8}\e \f yx\mod p.$
\endpro

\par Conjectures 2.7-2.9 have been checked for all primes $p<100,000$
and $q<100$.

 \pro{Conjecture 2.10 ([S3, Conjecture 9.5])} Let $p\e 1\mod 4$ be a prime,
$b\in\Bbb Z$, $2\nmid b$, $p\not=b^2+4$ and
$p=c^2+d^2=x^2+(b^2+4)y^2$ for some $c,d,x,y\in\Bbb Z$. Suppose $c\e
1\mod 4$ and all the odd parts of $d,x,y$ are numbers of the form
$4k+1$.
\par $(\t{\rm i})$ If $4\nmid xy$, then
$$\aligned\Ls{b+\f{cx}{dy}}2^{\f{p-1}4}&\e -\Ls{b-\f{cx}{dy}}2^{\f{p-1}4}
\\&\e\cases -(-1)^{\f d4}\f dc\mod p&\t{if $2\ \Vert\ x$
and $b\e 1,3\mod 8$,} \\(-1)^{\f d4}\f dc\mod p&\t{if $2\ \Vert\ x$
and $b\e 5,7\mod 8$,}
\\ 1\mod p&\t{if $2\ \Vert\ y$.}
\endcases\endaligned$$
\par $(\t{\rm ii})$ If $4\mid xy$, then
$$\aligned\Ls{b+\f{cx}{dy}}2^{\f{p-1}4}&\e \Ls{b-\f{cx}{dy}}2^{\f{p-1}4}
\\&\e \cases (-1)^{\f {d+y}4}
\mod p&\t{if $4\mid y$,}
\\ -(-1)^{\f x4}\f dc\mod p&\t{if $4\mid x$ and $b\e 1,3\mod 8$,}
\\ (-1)^{\f x4}\f dc\mod p&\t{if $4\mid x$ and $b\e 5,7\mod 8$.}
\endcases\endaligned$$
\endpro
\par We remark that Conjectures 2.1-2.10 have been solved by the author
under some restricted conditions in [S9].
 \subheading{3. Conjectures on Lucas
sequences}

 \pro{Conjecture 3.1 ([S3, Conjecture 9.4])} Let $p\e
1\mod 4$ be a prime, $b\in\Bbb Z$, $2\nmid b$ and
$p=c^2+d^2=x^2+(b^2+4)y^2\not=b^2+4$ for some $c,d,x,y\in\Bbb Z$.
Suppose $c\e 1\mod 4$ and all the odd parts of $d,x,y$ are numbers
of the form $4k+1$.
\par $(\t{\rm i})$ If $4\nmid xy$, then
$$U_{\f{p-1}4}(b,-1)\e \cases(-1)^{\f d4}\f{2y}x\mod p&\t{if $2\ \Vert\ x$
and $b\e 1,3\mod 8$,} \\-(-1)^{\f d4}\f{2y}x\mod p&\t{if $2\
\Vert\ x$ and $b\e 5,7\mod 8$,}
\\\f{2dy}{cx}\mod p&\t{if $2\ \Vert\ y$.}
\endcases$$
\par $(\t{\rm ii})$ If $4\mid xy$, then
$$V_{\f{p-1}4}(b,-1)\e \cases 2(-1)^{\f {d+y}4}
\mod p&\t{if $4\mid y$,}
\\ -2(-1)^{\f x4}\f dc\mod p&\t{if $4\mid x$ and $b\e 1,3\mod 8$.}
\\ 2(-1)^{\f x4}\f dc\mod p&\t{if $4\mid x$ and $b\e 5,7\mod 8$.}
\endcases$$
\endpro
\par Conjecture 3.1 has been checked for $1\le b<60$ and $p<20,000$.
When $p\e 1\mod 8$, $b=1,3$ and $4\mid y$, the conjecture
$V_{\f{p-1}4}(b,-1)\e 2(-1)^{\f{d+y}4}\mod p$ is equivalent to a
conjecture of E. Lehmer. See [L, Conjecture 4].
 \par By (1.3) and (1.4), Conjecture 3.1 is equivalent to Conjecture 2.11.
 By [S3], Conjectures 2.3 and 2.5 are consequences of Conjecture
 3.1.

  \pro{Conjecture 3.2 ([S3, Conjecture 9.11])} Let $p\e 1\mod 4$ be a
prime, $b\in\Bbb Z$, $b\e 4\mod{8}$, $p\not=b^2/4+1$ and
$p=c^2+d^2=x^2+(1+b^2/4)y^2$ for some $c,d,x,y\in\Bbb Z$. Suppose
$c\e 1\mod 4$ and all the odd parts of $d,x,y$ are numbers of the
form $4k+1$. Then
$$U_{\f{p-1}4}(b,-1)\e\cases (-1)^{\f{b+4}8+\f d4}\f yx\mod p&\t{if $2\
\Vert\ x$,}\\\f{dy}{cx}\mod p&\t{if $2\ \Vert\ y$,}
\\0\mod p&\t{if $4\mid xy$}\endcases$$ and
$$V_{\f{p-1}4}(b,-1)\e\cases 2(-1)^{\f d4+\f y4}\mod p&\t{if $4\mid y$,}
\\2(-1)^{\f{b-4}8+\f x4}\f dc\mod p&\t{if $4\mid x$,}
\\0\mod p&\t{if $4\nmid xy$.}
\endcases$$
\endpro
\par Conjecture 3.2 has been checked for $1\le b\le 100$ and $p<20,000$.

\pro{Conjecture 3.3 ([S3, Conjecture 9.14])} Let $p\e 1\mod 4$ be a
prime, $b\in\Bbb Z$, $8\mid b$, $p\not=b^2/4+1$ and
$p=c^2+d^2=x^2+(1+b^2/4)y^2$ for some $c,d,x,y\in\Bbb Z$. Suppose
$c\e 1\mod 4$ and all the odd parts of $d,x,y$ are numbers of the
form $4k+1$. Then
$$U_{\f{p-1}4}(b,-1)\e \cases 0\mod p&\t{if $4\mid xy$,}
\\-(-1)^{(\f b8-1)y}\f{dy}{cx}\mod p&\t{if $4\nmid xy$}\endcases$$
and
$$V_{\f{p-1}4}(b,-1)\e \cases 2(-1)^{\f d4+\f {xy}4+\f b8y}\mod p
&\t{if $4\mid xy$,}\\0\mod p&\t{if $4\nmid xy$.}
\endcases$$\endpro
\par Conjecture 3.3 has been checked for $1\le b<100$ and $p<20,000$. By
[S3], Conjecture 3.3 implies Conjecture 2.6.

\pro{Conjecture 3.4 ([S3, Conjecture 9.17])} Let $p\e 1\mod 4$ be a
prime, $b\in\Bbb Z$, $b\e 2\mod{4}$, $p\not=b^2/4+1$ and
$p=c^2+d^2=x^2+(1+b^2/4)y^2$ for some $c,d,x,y\in\Bbb Z$. Suppose
$c\e 1\mod 4$, $x=2^{\a}x_0$, $y=2^{\beta}y_0$ and $x_0\e y_0\e
1\mod 4$. Then
$$U_{\f{p-1}4}(b,-1)\e\cases (-1)^{\f{b-2}4+\f d4}\f yx\mod p&\t{if $2\
\Vert\ y$,}\\0\mod p&\t{if $4\mid y$}\endcases$$ and
$$V_{\f{p-1}4}(b,-1)\e\cases 0\mod p&\t{if $2\ \Vert\ y$,}
\\2(-1)^{\f d4+\f y4}\mod p&\t{if $4\mid y$.}
\endcases$$
\endpro
\par Conjecture 3.4 has been checked for $1\le b<100$ and $p<20,000$. By
[S3], Conjecture 3.4 implies Conjecture 2.4.

\pro{Conjecture 3.5 ([S5, Conjecture 4.1])} Let $p$ be an odd prime
and $k\in\Bbb N$ with $2\nmid k$. Suppose $p=x^2+(k^2+1)y^2$ for
some $x,y\in\Bbb Z$. Then
$$V_{\f{p+1}4}(2k,-1)\e \cases
(-1)^{\f{(\f{p-1}2y)^2-1}8}(-2)^{\f{p+1}4}\mod p&\t{if $k\e 5,7\mod
8$,}
\\-(-1)^{\f{(\f{p-1}2y)^2-1}8}(-2)^{\f{p+1}4}\mod p&\t{if $k\e 1,3\mod 8$.}
\endcases$$\endpro
\par Conjecture 3.5 has been checked for all $k<60$ and
$p<20,000$. The case $k=1$ has been solved by the author in [S3].

 \pro{Conjecture 3.6 ([S5, Conjecture 4.2])} Let $p$ be an odd prime and $k\in\Bbb N$
with $2\nmid k$. Suppose $2p=x^2+(k^2+4)y^2$ for some $x,y\in\Bbb
Z$.
\par $(\t{\rm i})$ If $k\e 1,3\mod 8$, then
$$V_{\f{p+1}4}(k,-1)\e \cases
(-1)^{\f{(\f{p-1}2y)^2-1}8}(-2)^{\f{p+1}4}\mod p&\t{if $k\e 1,11\mod
{16}$,}
\\-(-1)^{\f{(\f{p-1}2y)^2-1}8}(-2)^{\f{p+1}4}\mod p&\t{if $k\e 3,9\mod {16}$.}
\endcases$$
\par $(\t{\rm ii})$ If $k\e 5,7\mod 8$, then
$$V_{\f{p+1}4}(k,-1)\e \cases
(-1)^{\f{(\f{p-1}2y)^2-1}8}2^{\f{p+1}4}\mod p&\t{if $k\e 5,15\mod
{16}$,}
\\-(-1)^{\f{(\f{p-1}2y)^2-1}8}2^{\f{p+1}4}\mod p&\t{if $k\e 7,13\mod {16}$.}
\endcases$$\endpro
\par Conjecture 3.6 has been checked for all $k<60$ and
$p<20,000$. The case $k=1$ was posed by the author in [S1] and
solved by C.N. Beli in [B].

 \pro{Conjecture
3.7 ([S6, Conjecture 5.1])} Let $p>7$ be a prime such that $p\e
1,2,4\mod 7$ and so $p=C^2+7D^2$ with $C,D\in\Bbb Z$.
\par $(\t{\rm i})$ If $p\e 1\mod 4$ and $C\e 1\mod 4$, then
$$\align&\sum_{k=0}^{(p-1)/2}\binom{2k}k^2U_k(16,1)
\e 0\mod{p^2},
\\&\sum_{k=0}^{(p-1)/2}\binom{2k}k^2V_k(16,1) \e (-1)^{\f{p-1}4}(4C-\f
pC)\mod{p^2}.\endalign$$
\par $(\t{\rm ii})$ If $p\e 3\mod 4$ and $D\e 1\mod 4$, then
$$\align&\sum_{k=0}^{(p-1)/2}\binom{2k}k^2U_k(16,1)\e
(-1)^{\f{p-3}4}\Big(\f{16}3D-\f {4p}{21D}\Big)\mod{p^2},
\\&\sum_{k=0}^{(p-1)/2}\binom{2k}k^2V_k(16,1)
\e (-1)^{\f{p+1}4}\Big(84D-\f {3p}D\Big)\mod{p^2}.\endalign$$
\endpro
\par In [S6], the author proved the congruences modulo $p$.
\pro{Conjecture 3.8 ([S6, Conjecture 5.2])} Let $p>3$ be a prime
such that $p\e 1,3,4,9,10,12\mod{13}$.  Then
$$\sum_{k=0}^{(p-1)/2}\binom{2k}k^2U_k(11,1)\e 0\mod p.$$
\endpro

\subheading{4. Conjectures on supercongruences}
\par In 2003, Rodriguez-Villegas[RV] posed many conjectures on supercongruences. In particular,
he conjectured that for any prime $p>3$,
$$\aligned&\sum_{k=0}^{p-1}\f{(6k)!}{1728^k(3k)!k!^3}\e \cases
\sls p3(4x^2-2p)\mod{p^2}&\t{if $p=x^2+y^2$ with $2\nmid x$,}
\\0\mod {p^2}&\t{if $p\e 3\mod 4$}\endcases
\\&\sum_{k=0}^{p-1}\f{(4k)!}{256^kk!^4}\\&\e\cases 4x^2-2p\mod {p^2}
&\t{if $p\e 1,3\mod 8$ and so $p=x^2+2y^2$,}
\\ 0\mod{p^2}&\t{if $p\e 5,7\mod 8$}
\endcases\endaligned$$
and
$$\sum_{k=0}^{p-1}\f{\binom{2k}k^2\binom{3k}k}{108^k}
\e\cases 4x^2-2p\mod {p^2}&\t{if $3\mid p-1$ and $p=x^2+3y^2$,}
\\0\mod {p^2}&\t{if $3\mid p-2$.}
\endcases$$
The three conjectures have been solved by E. Mortenson[M] and
Zhi-Wei Sun[Su2]. Recently my twin brother Zhi-Wei Sun made a lot of
conjectures on supercongruences. Inspired by his work in [Su1], the
author made the following conjectures.
 \pro{Conjecture 4.1 ([S4, 2.1])}
Let $p\e 1\mod 4$ be a prime and so $p=x^2+y^2$ with $2\nmid x$.
Then
$$\sum_{k=0}^{p-1}\f{(4k)!}{648^kk!^4}\e 4x^2-2p\mod {p^2}.$$
\endpro
\par In [S6], the author proved the congruence modulo $p$. When
$p\e 3\mod 4$, the author proved
$\sum_{k=0}^{p-1}\f{(4k)!}{648^kk!^4}\e 0\mod{p^2}$ in [S8].

\pro{Conjecture 4.2 ([S4, 2.2])} Let $p\e 1\mod 3$ be a prime and so
$p=x^2+3y^2$. Then
$$\sum_{k=0}^{p-1}\f{(4k)!}{(-144)^kk!^4}\e 4x^2-2p\mod {p^2}.$$
\endpro
\par In [S6], the author proved the congruence modulo $p$. When
$p\e 5\mod 6$, the author proved
$\sum_{k=0}^{p-1}\f{(4k)!}{(-144)^kk!^4}\e 0\mod{p^2}$ in [S8].
\pro{Conjecture 4.3 ([S4, 2.3])} Let $p\e 1,2,4\mod 7$ be an odd
prime and so $p=x^2+7y^2$. Then
$$\sum_{k=0}^{p-1}\f{(4k)!}{(-3969)^kk!^4}\e 4x^2-2p\mod {p^2}.$$
\endpro
\par In [S6], the author proved the congruence modulo $p$. When
$p\e 3,5,6\mod 7$, the author proved
$\sum_{k=0}^{p-1}\f{(4k)!}{(-3969)^kk!^4}\e 0\mod{p^2}$ in [S8].
  \pro{Conjecture 4.4 ([S4, 2.4])} Let $p\e 1\mod 4$ be a
prime and so $p=x^2+y^2$ with $2\nmid x$. Then
$$\sum_{k=0}^{p-1}\f{(6k)!}{66^{3k}(3k)!k!^3}
\e \Ls p{33}(4x^2-2p)\mod {p^2}.
$$\endpro
\par In [S7], the author proved the congruence modulo $p$. When
$p\e 3\mod 4$, the author proved
$\sum_{k=0}^{p-1}\f{(6k)!}{66^{3k}(3k)!k!^3}\e 0\mod{p^2}$ in [S8].

\pro{Conjecture 4.5 ([S4, 2.5])} Let $p\e 1,3\mod 8$ be a prime and
so $p=x^2+2y^2$. Then
$$\sum_{k=0}^{p-1}\f{(6k)!}{20^{3k}(3k)!k!^3}\e
\Ls {-5}p(4x^2-2p)\mod {p^2}.$$ \endpro
\par In [S7], the author proved the congruence modulo $p$. When
$p\e 5,7\mod 8$, the author proved
$\sum_{k=0}^{p-1}\f{(6k)!}{20^{3k}(3k)!k!^3}\e 0\mod{p^2}$ in [S8].

\pro{Conjecture 4.6 ([S4, 2.6])} Let $p\e 1\mod 3$ be a prime and so
$p=x^2+3y^2$. Then
$$\sum_{k=0}^{p-1}\f{(6k)!}{54000^k(3k)!k!^3}\e
\Ls p5(4x^2-2p)\mod {p^2}.$$ \endpro
\par In [S7], the author proved the congruence modulo $p$. When
$p\e 5\mod 6$, the author proved
$\sum_{k=0}^{p-1}\f{(6k)!}{54000^{3k}(3k)!k!^3}\e 0\mod{p^2}$ in
[S8].

\pro{Conjecture 4.7 ([S4, 2.7])} Let $p>5$ be a prime. Then
$$\aligned&\sum_{k=0}^{p-1}\f{(6k)!}{(-12288000)^k(3k)!k!^3}\\&\e\cases \sls {10}p(L^2-2p)\mod {p^2}
&\t{if $p\e 1\mod 3$ and so $4p=L^2+27M^2$,}
\\ 0\mod{p^2}&\t{if $p\e 2\mod 3$.}
\endcases\endaligned$$\endpro

\pro{Conjecture 4.8 ([S4, 2.8])} Let $p\e 1,2,4\mod 7$ be an odd
prime and so $p=x^2+7y^2$. Then
$$\sum_{k=0}^{p-1}\f{(6k)!}{(-15)^{3k}(3k)!k!^3}
\e \Ls p{15}(4x^2-2p)\mod {p^2}.$$ \endpro
\par In [S7], the author proved the congruence modulo $p$. When
$p\e 3,5,6\mod 7$, the author proved
$\sum_{k=0}^{p-1}\f{(6k)!}{(-15)^{3k}(3k)!k!^3}\e 0\mod{p^2}$ in
[S8].

\pro{Conjecture 4.9 ([S4, 2.9])} Let $p\e 1,2,4\mod 7$ be an odd
prime and so $p=x^2+7y^2$. Then
$$\sum_{k=0}^{p-1}\f{(6k)!}{255^{3k}(3k)!k!^3}\e
\Ls p{255}(4x^2-2p)\mod {p^2}.$$\endpro
 \par In [S7], the author proved the congruence modulo $p$. When
$p\e 3,5,6\mod 7$, the author proved
$\sum_{k=0}^{p-1}\f{(6k)!}{255^{3k}(3k)!k!^3}\e 0\mod{p^2}$ in [S8].

\pro{Conjecture 4.10 ([S4, 2.10])} Let $p>3$ be a prime. Then
$$\aligned&\sum_{k=0}^{p-1}\f{\binom{2k}k^2\binom{3k}k}{1458^k}\\&\e\cases 4x^2-2p\mod {p^2}
&\t{if $p\e 1\mod 3$ and so $p=x^2+3y^2$,}
\\ 0\mod{p^2}&\t{if $p\e 2\mod 3$.}
\endcases\endaligned$$\endpro

\pro{Conjecture 4.11 ([S4, 2.11])} Let $p>5$ be a prime. Then
$$\aligned&\sum_{k=0}^{p-1}\f{\binom{2k}k^2\binom{3k}k}{15^{3k}}\\&\e\cases 4x^2-2p\mod {p^2}
&\t{if $p\e 1,4\mod {15}$ and so $p=x^2+15y^2$,}
\\2p-12x^2\mod {p^2}
&\t{if $p\e 2,8\mod {15}$ and so $p=3x^2+5y^2$,}
\\ 0\mod{p^2}&\t{if $p\e 7,11,13,14\mod {15}$.}
\endcases\endaligned$$\endpro

\pro{Conjecture 4.12 ([S4, 2.12])} Let $p>3$ be a prime. Then
$$\aligned&\sum_{k=0}^{p-1}\f{\binom{2k}k^2\binom{3k}k}{(-8640)^k}\\&\e
\cases 4x^2-2p\mod{p^2}&\t{if $3\mid p-1$, $p=x^2+3y^2$ and
$10^{\f{p-1}3}\e 1\mod p$,}
\\p-2x^2\pm 6xy\mod{p^2}&\t{if $3\mid p-1$, $p=x^2+3y^2$ and $10^{\f{p-1}3}\e \f
12(-1\mp \f xy)\mod p$,}
\\0\mod{p^2}&\t{if $3\mid p-2$.}\endcases\endaligned$$
\endpro
\par If $p$ is a prime of the form $3k+1$ and so $p=x^2+3y^2$ with
$x,y\in\Bbb Z$, by [S2, Corollary 4.4] we have $10^{\f{p-1}3}\e
1\mod p$ if and only if $5\mid xy$.

\pro{Conjecture 4.13 ([S6, Conjecture 5.3])} Let $p\e 1\mod{12}$ be
a prime and $p=a^2+4b^2=A^2+3B^2$ with $a\e A
 \e 1\mod 4$. Then
 $$\sum_{k=0}^{(p-1)/6}\f{\binom{6k}{3k}^2}{(-16)^{3k}}
 \e (-1)^{\f{p-1}4}\f 13\Big(2a+4A-\f p{2a}-\f pA\Big)\mod {p^2}.$$
 \endpro
\par In [S6], the author proved the congruence modulo $p$.

\par For an integer $m$ and odd prime $p$ with $p\nmid m$ let

$$S_p(m)=\sum_{n=0}^{p-1}\sum_{k=0}^n\b nk^4m^k.$$
Recently Zhi-Wei Sun [Su5] investigated congruences $S_p(m)\mod
{p^2}$ and revealed the connections with binary quadratic forms and
series for $\f 1{\pi}$. Now we introduce the sum
$$Z_p(m)=\sum_{n=0}^{p-1}\f{\b{2n}n}{m^n}\sum_{k=0}^n\b nk^3.$$
Then we have the following conjectures concerning $Z_p(m)$ modulo
$p^2$.
 \pro{Conjecture 4.14} Let $p$ be an odd prime. Then
$$Z_p(-16)\e \cases 4x^2-2p\mod {p^2}&\t{if $p=x^2+y^2\e 1\mod{12}$ with $6\mid y$,}
\\2p-4x^2\mod{p^2}&\t{if $p=x^2+y^2\e 1\mod{12}$ with $6\mid x-3$,}
\\4\sls {xy}3xy\mod{p^2}&\t{if $p=x^2+y^2\e 5\mod{12}$,}
\\0\mod{p^2}&\t{if $p\e 3\mod 4$.}\endcases
$$\endpro
 \pro{Conjecture 4.15} Let $p$ be an odd prime. Then
$$Z_p(96)\e \cases \sls p3(4x^2-2p)\mod{p^2}&\t{if $p=x^2+2y^2\e
1,3\mod 8$,}\\0\mod{p^2}&\t{if $p\e 5,7\mod 8$.}
\endcases$$
\endpro
\pro{Conjecture 4.16} Let $p>5$ be a prime. Then
$$Z_p(50)\e \cases 4x^2-2p\mod{p^2}&\t{if $p=x^2+3y^2\e
1\mod 3$,}\\0\mod{p^2}&\t{if $p\e 2\mod 3$.}
\endcases$$
\endpro
\pro{Conjecture 4.17} Let $p>5$ be a prime. Then
$$Z_p(16)\e \cases 4x^2-2p\mod{p^2}&\t{if $p\e
1,9\mod{20}$ and so $p=x^2+5y^2$,}\\2x^2-2p\mod{p^2}&\t{if $p\e
3,7\mod{20}$ and so $2p=x^2+5y^2$,}
\\0\mod{p^2}&\t{if $p\e 11,13,17,19\mod{20}$.}\endcases$$\endpro

\pro{Conjecture 4.18} Let $p>3$ be a prime. Then
$$Z_p(32)\e \cases 4x^2-2p\mod{p^2}&\t{if $p\e
1,7\mod{24}$ and so $p=x^2+6y^2$,}\\8x^2-2p\mod{p^2}&\t{if $p\e
5,11\mod{24}$ and so $p=2x^2+3y^2$,}\\0\mod{p^2}&\t{if $p\e
13,17,19,23\mod{24}$.}\endcases$$\endpro

\pro{Conjecture 4.19} Let $p>7$ be a prime. Then
$$Z_p(5)\e Z_p(-49)\e \cases 4x^2-2p\mod{p^2}&\t{if $p=x^2+15y^2\e
1,19\mod{30}$,}\\2p-12x^2\mod{p^2}&\t{if $p=3x^2+5y^2\e
17,23\mod{30}$,}
\\0\mod{p^2}&\t{if $p\e 7,11,13,29\mod{30}$.}\endcases$$\endpro

\pro{Conjecture 4.20} Let $b\in\{7,11,19,31,59\}$ and
$$f(b)=\cases -112&\t{if $b=7$,}
\\-400&\t{if $b=11$,}
\\-2704&\t{if $b=19$,}
\\-24304&\t{if $b=31$,}
\\-1123600&\t{if $b=59$.}
\endcases$$
If $p$ is a prime with $p\not= 2,3,b$ and $p\nmid f(b)$, then
$$Z_p(f(b))\e \cases 4x^2-2p\mod{p^2}&\t{if $p=x^2+3by^2$,}
\\2p-12x^2\mod{p^2}&\t{if $p=3x^2+by^2$,}
\\2x^2-2p\mod{p^2}&\t{if $2p=x^2+3by^2$,}
\\2p-6x^2\mod{p^2}&\t{if $2p=3x^2+by^2$,}
\\0\mod{p^2}&\t{if $\sls{-3b}p=-1$.}
\endcases$$\endpro

\pro{Conjecture 4.21} Let $b\in\{5,7,13,17\}$ and
$$f(b)=\cases 320&\t{if $b=5$,}
\\896&\t{if $b=7$,}
\\10400&\t{if $b=13$,}
\\39200&\t{if $b=17$.}
\endcases$$
If $p$ is a prime with $p\not= 2,3,b$ and $p\nmid f(b)$, then
$$Z_p(f(b))\e \cases 4x^2-2p\mod{p^2}&\t{if $p=x^2+6by^2$,}
\\8x^2-2p\mod{p^2}&\t{if $p=2x^2+3by^2$,}
\\2p-12x^2\mod{p^2}&\t{if $p=3x^2+2by^2$,}
\\2p-24x^2\mod{p^2}&\t{if $p=6x^2+by^2$,}
\\0\mod{p^2}&\t{if $\sls{-6b}p=-1$.}
\endcases$$\endpro

\pro{Conjecture 4.22} Let $p>3$ be a prime. Then
$$\align&\sum_{n=0}^{p-1}\f{9n+4}{5^n}\b{2n}n\sum_{k=0}^n\b nk^3
\e 4\Ls p5p\mod{p^2}\qtq{for}p>5,
\\&\sum_{n=0}^{p-1}\f{5n+2}{16^n}\b{2n}n\sum_{k=0}^n\b nk^3
\e 2p\mod{p^2},
\\&\sum_{n=0}^{p-1}\f{9n+2}{50^n}\b{2n}n\sum_{k=0}^n\b nk^3
\e 2\Ls {-1}pp\mod{p^2}\qtq{for}p\not=5,
\\&\sum_{n=0}^{p-1}\f{5n+1}{96^n}\b{2n}n\sum_{k=0}^n\b nk^3
\e \Ls {-2}pp\mod{p^2},
\\&\sum_{n=0}^{p-1}\f{6n+1}{320^n}\b{2n}n\sum_{k=0}^n\b nk^3
\e \Ls p{15}p\mod{p^2}\qtq{for}p\not=5,
\\&\sum_{n=0}^{p-1}\f{90n+13}{896^n}\b{2n}n\sum_{k=0}^n\b nk^3
\e 13\Ls p{7}p\mod{p^2}\qtq{for}p\not=7,
\\&\sum_{n=0}^{p-1}\f{102n+11}{10400^n}\b{2n}n\sum_{k=0}^n\b nk^3
\e 11\Ls p{39}p\mod{p^2}\qtq{for}p\not=5,13
\endalign$$
\endpro

\pro{Conjecture 4.23} Let $p>3$ be a prime. Then
$$\align
\\&\sum_{n=0}^{p-1}\f{3n+1}{(-16)^n}\b{2n}n\sum_{k=0}^n\b nk^3
\e \Ls {-1}pp\mod{p^2},
\\&\sum_{n=0}^{p-1}\f{15n+4}{(-49)^n}\b{2n}n\sum_{k=0}^n\b nk^3
\e 4\Ls p3p\mod{p^2}\qtq{for}p\not=5,7,
\\&\sum_{n=0}^{p-1}\f{9n+2}{(-112)^n}\b{2n}n\sum_{k=0}^n\b nk^3
\e 2\Ls p7p\mod{p^2}\qtq{for}p\not=7,
\\&\sum_{n=0}^{p-1}\f{99n+17}{(-400)^n}\b{2n}n\sum_{k=0}^n\b nk^3
\e 17\Ls{-1}pp\mod{p^2},
\\&\sum_{n=0}^{p-1}\f{855n+109}{(-2704)^n}\b{2n}n\sum_{k=0}^n\b nk^3
\e 109p\Ls{-1}p\mod{p^2}\qtq{for}p\not=13,
\\&\sum_{n=0}^{p-1}\f{585n+58}{(-24304)^n}\b{2n}n\sum_{k=0}^n\b nk^3
\e 58p\Ls{-31}p\mod{p^2}\qtq{for}p\not=7,31.
\endalign$$
\endpro
\pro{Conjecture 4.24} Let $p>5$ be a prime. Then
$$\align &\sum_{k=0}^{p-1}\f{63k+8}{(-15)^{3k}}
\b{2k}k\b{3k}k\b{6k}{3k}\e 8p\Ls {-15}p\mod{p^2},
\\&\sum_{k=0}^{p-1}\f{133k+8}{255^{3k}}
\b{2k}k\b{3k}k\b{6k}{3k}\e 8p\Ls {-255}p\mod{p^2}\qtq{for}p\not=17,
\\&\sum_{k=0}^{p-1}\f{28k+3}{20^{3k}}
\b{2k}k\b{3k}k\b{6k}{3k}\e 3p\Ls{-5}p\mod{p^2},
\\&\sum_{k=0}^{p-1}\f{63k+5}{66^{3k}}
\b{2k}k\b{3k}k\b{6k}{3k}\e 5p\Ls{-33}p\mod{p^2}\qtq{for}p\not=11,
\\&\sum_{k=0}^{p-1}\f{11k+1}{54000^k}
\b{2k}k\b{3k}k\b{6k}{3k}\e p\Ls {-15}p\mod{p^2},
\\&\sum_{k=0}^{p-1}\f{506k+31}{(-12288000)^k}
\b{2k}k\b{3k}k\b{6k}{3k}\e 31p\Ls{-30}p\mod{p^2}.
\endalign$$
\endpro

\pro{Conjecture 4.25} Let $p>5$ be a prime. Then
$$\aligned\sum_{k=0}^{p-1}\b{-\f 13}k\b{-\f 16}k(-4)^k
&\e \Ls 5p5^{\f{1-\sls p3}2}\sum_{k=0}^{p-1}\b{-\f 23}k\b{-\f
56}k(-4)^k\\& \e\cases \sls x3(2x-\f p{2x})\mod {p^2}&\t{if
$p=x^2+15y^2$,}
\\ -\sls x3(10x-\f p{2x})\mod {p^2}&\t{if $p=5x^2+3y^2$,}
\\0\mod p&\t{if $p\e 17,23\mod{30}$}
\endcases\endaligned$$ and so
$$2x\Ls x3\e
\sum_{k=0}^{(p-5)/6}\b{\f{p-2}3}k\b{\f{p-5}6}k(-4)^k\mod
p\qtq{for}p=5x^2+3y^2.$$
\endpro

\pro{Conjecture 4.26} Let $p>3$ be a prime. Then
$$\aligned&\sum_{k=0}^{p-1}\f{\b{-\f 13}k\b{-\f 16}k}{2^k}
\e \Ls 2p2^{\f{\sls p3-1}2}\sum_{k=0}^{p-1}\f{\b{-\f 23}k\b{-\f
56}k}{2^k}
\\&\e\cases \sls x3(2x-\f p{2x})\mod {p^2}&\t{if $p=x^2+6y^2\e 1,7\mod
{24}$,}\\\sls x3(2x-\f p{4x})\mod {p^2}&\t{if $p=2x^2+3y^2\e
5,11\mod {24}$,}
\\0\mod p&\t{if $p\e 13,19\mod {24}$}\endcases\endaligned$$
and so
$$x\Ls x3\e -\f 14
\sum_{k=0}^{(p-5)/6}\b{\f{p-2}3}k\b{\f{p-5}6}k\f 1{2^k}\mod
p\qtq{for}p=2x^2+3y^2.$$
\endpro

\pro{Conjecture 4.27} Let $p>3$ be a prime. Then
$$\aligned\sum_{k=0}^{p-1}\f{\b{-\f 13}k\b{-\f 16}k}{(-16)^k}
&\e\Ls{17}p\Ls{17}{16}^{\f{1-\sls p3}2}\sum_{k=0}^{p-1}\f{\b{-\f
23}k\b{-\f 56}k}{(-16)^k} \\&\e\cases -\sls x3(x-\f px)\mod
{p^2}&\t{if $4p=x^2+51y^2$,}
\\ \f 14\sls x3(17x-\f px)\mod{p^2}&\t{if $4p=17x^2+3y^2$,}
\\0\mod p&\t{if $\sls p3=-\sls p{17}=1$}\endcases\endaligned$$
and so
$$x\Ls x3\e -\f 14
\sum_{k=0}^{(p-5)/6}\b{\f{p-2}3}k\b{\f{p-5}6}k\f 1{(-16)^k}\mod
p\qtq{for}4p=17x^2+3y^2.$$

\endpro

\pro{Conjecture 4.28} Let $p>3$ be a prime. Then
$$\aligned\sum_{k=0}^{p-1}\f{\b{-\f 13}k\b{-\f 16}k}{(-1024)^k}
&\e \Ls{41}p\Ls{1025}{1024}^{\f{1-\sls p3}2}
\sum_{k=0}^{p-1}\f{\b{-\f 23}k\b{-\f 56}k}{(-1024)^k}
 \\&\e\cases -\sls x3(x-\f px)\mod {p^2}&\t{if $4p=x^2+123y^2$,}
\\ \f 5{32}\sls x3(41x-\f px)\mod{p^2}&\t{if $4p=41x^2+3y^2$,}
\\0\mod p&\t{if $\sls p3=-\sls p{41}=1$}\endcases\endaligned$$

and so
$$x\Ls x3\e
-\f 5{32}\sum_{k=0}^{(p-5)/6}\b{\f{p-2}3}k\b{\f{p-5}6}k\f
1{(-1024)^k}\mod p\qtq{for}4p=41x^2+3y^2.$$

\endpro

\pro{Conjecture 4.29} Let $p>3$ be a prime. Then
$$\aligned\sum_{k=0}^{p-1}\f{\b{-\f 13}k\b{-\f 16}k}{(-250000)^k}
&\e \Ls {89}p\Ls{250001}{250000}^{\f{1-\sls
p3}2}\sum_{k=0}^{p-1}\f{\b{-\f 23}k\b{-\f 56}k}{(-250000)^k}
 \\&\e\cases -\sls x3(x-\f px)\mod {p^2}&\t{if $4p=x^2+267y^2$,}
\\ \f{53}{500}\sls x3(89x-\f px)\mod{p^2}&\t{if $4p=89x^2+3y^2$,}
\\0\mod p&\t{if $\sls p3=-\sls p{89}=1$}\endcases\endaligned$$
and so
$$x\Ls x3\e -\f{53}{500}
\sum_{k=0}^{(p-5)/6}\b{\f{p-2}3}k\b{\f{p-5}6}k\f 1{(-250000)^k}\mod
p\ \t{for}\ 4p=89x^2+3y^2.$$
\endpro
\pro{Conjecture 4.30}  Let $p>5$ be a prime. Then
$$\aligned&\sum_{k=0}^{p-1}\b{-\f 13}k\b{-\f 16}k\f 1{(-80)^k}
\\&\e \Ls 5p
\sum_{k=0}^{p-1}\b{-\f 23}k\b{-\f 56}k\f 1{(-80)^k}
\\&\e \cases x-\f px\mod
{p^2}&\t{if $p\e 1,19\mod {30}$ and so $4p=x^2+75y^2$ with $3\mid
x-2$,}
\\5x-\f p{5x}\mod
{p^2}&\t{if $p\e 7,13\mod {30}$ and so $4p=25x^2+3y^2$ with $3\mid
x-1$,}
 \\0\mod p&\t{if
$p\e 17,23\mod{30}$}\endcases\endaligned$$
\endpro

\pro{Conjecture 4.31} Let $p>11$ be a prime such that $\sls p{11}=1$
and so $4p=x^2+11y^2$. Then

$$\aligned\sum_{k=0}^{p-1}\b{-\f 13}k\b{-\f 16}k\Ls {27}{16}^k
&\e\Big(-\f{11}{16}\Big)^{\f{1-\sls p3}2} \sum_{k=0}^{p-1}\b{-\f
23}k\b{-\f 56}k\Ls {27}{16}^k
\\& \e \cases -\sls{-11+x/y}p\sls x{11}(x-\f px)\mod {p^2}&\t{if
$3\mid p-1$,}
\\-\f 14\sls{-11+\sls x{11}x/y}p(11y-\f py)\mod
{p^2}&\t{if $3\mid p-2$}
 \endcases\endaligned$$
and so
$$\align &y\Ls{-11+\sls x{11}x/y}p\\&\e \f 14\sum_{k=0}^{(p-5)/6}
\b{\f{p-2}3}k\b{\f{p-5}6}k\Ls{27}{16}^k\mod p\qtq{for}4p=x^2+11y^2\e
2\mod 3.\endalign$$

\endpro

\pro{Conjecture 4.32} Let $p\e 1\mod 3$ be a prime and so
$4p=L^2+27M^2$ with $3\mid L-2$. Then
$$\sum_{k=0}^{p-1}\b{-\f 13}k\b{-\f 16}k\Big(-\f 9{16}\Big)^k
\e \sum_{k=0}^{p-1}\b{-\f 23}k\b{-\f 56}k\Big(-\f 9{16}\Big)^k\e
L-\f pL\mod {p^2}.$$\endpro

\pro{Conjecture 4.33} Let $p\e 1,3\mod 8$ be a prime and so
$p=c^2+2d^2$ with $4\mid c-1$. Then
$$\aligned&\sum_{k=0}^{p-1}\b{-\f 13}k\b{-\f 16}k\Big(\f {27}2\Big)^k
\\&\e \Big(-\f{25}2\Big)^{\f{1-\sls p3}2}\sum_{k=0}^{p-1}\b{-\f
23}k\b{-\f 56}k\Ls {27}2^k
\\&\e\cases (-1)^{[\f{p}8]}\sls{-2-c/d}p(2c-\f p{2c}) \mod {p^2}&\t{if
$p\e 1,19\mod{24}$,}
\\(-1)^{[\f{p}8]}\sls{2+c/d}p(10d-\f{5p}{4d}) \mod {p^2}
&\t{if $p\e 11,17\mod{24}$.}
\endcases\endaligned$$\endpro

\pro{Conjecture 4.34} Let $p$ be a prime such that $p\e 1\mod 3$,
  $4p=L^2+27M^2(L,M\in\Bbb Z)$ and $L\e 2\mod 3$.  Then
$$\sum_{k=0}^{p-1}\b{-\f 13}k^29^k
\e L-\f pL\mod{p^2}$$\endpro

\pro{Conjecture 4.35} Let $p$ be a prime of the form $4k+1$  and so
$p=x^2+y^2$ with $4\mid x-1$. Then
$$\sum_{k=0}^{p-1}\b{-\f 14}k^2(-8)^k
\e(-1)^{\f{p-1}4}(2x-\f p{2x})\mod {p^2}$$ and
$$\aligned\sum_{k=0}^{p-1}\f{\b{-\f 14}k^2}{(-8)^k}
\e\cases (-1)^{\f y4}(2x-\f p{2x})\mod {p^2} &\t{if $p\e 1\mod
8$,}\\(-1)^{\f{y-2}4} (2y-\f p{2y})\mod {p^2}&\t{if $p\e 5\mod 8$.}
\endcases\endaligned$$\endpro

\pro{Conjecture 4.36} Let $p\e 1\mod 3$ be a prime and so
$p=A^2+3B^2$. Then
$$\aligned&\sum_{k=0}^{p-1}\b{-\f 14}k^24^k\e
\sum_{k=0}^{p-1}\b{-\f 14}k\b{-\f 12}k(-8)^k
\\&\e\cases (-1)^{\f{p-1}4+\f{A-1}2}
(2A-\f p{2A})\mod {p^2}&\t{if $p\e 1\mod {12}$,}
\\(-1)^{\f{p+1}4+\f{B-1}2}
(6B-\f p{2B})\mod {p^2}&\t{if $p\e 7\mod {12}$.}
\endcases\endaligned$$ and
$$\aligned\sum_{k=0}^{p-1}\f{\b{-\f 14}k^2}{4^k}
\e\cases (-1)^{\f{A-1}2}(2A-\f p{2A})\mod {p^2}&\t{if $p\e
1\mod{12}$,}
\\(-1)^{\f{B-1}2}(3B-\f p{4B})\mod {p^2}&\t{if $p\e
7\mod{12}$.}
\endcases\endaligned$$\endpro

\pro{Conjecture 4.37} Let $p>2$ be a prime such that $p\e 1,2,4\mod
7$ and so $p=x^2+7y^2$. Then
$$\aligned\sum_{k=0}^{p-1}\b{-\f 14}k^264^k
\e\cases \sls 2p(-1)^{\f{x-1}2}(2x-\f p{2x})\mod {p^2}&\t{if $p\e
1\mod 4$,}\\\sls 2p(-1)^{\f{y-1}2}(42y-\f {3p}{2y})\mod {p^2}&\t{if
$p\e 3\mod 4$}
\endcases\endaligned$$
and
$$\aligned\sum_{k=0}^{p-1}\f{\b{-\f 14}k^2}{64^k}
\e\cases (-1)^{\f{x-1}2}(2x-\f p{2x})\mod {p^2}&\t{if $p\e 1\mod
4$,}\\\f{3}4(-1)^{\f{y-1}2}(7y-\f p{4y})\mod {p^2}&\t{if $p\e 3\mod
4$.}
\endcases\endaligned$$\endpro

\pro{Conjecture 4.38} Let $p$ be a prime such that $p\e 5,7\mod 8$.
Then
$$\sum_{k=0}^{p-1}(-1)^k\b{-\f 14}k^2\e
\cases (-1)^{\f{x+1}2}(2x-\f p{2x})\mod {p^2} &\t{if $p=x^2+2y^2\e
1\mod 8$,}\\(-1)^{\f{y-1}2}(4y-\f p{2y})\mod {p^2} &\t{if
$p=x^2+2y^2\e 3\mod 8$},\\0\mod p&\t{if $p\e 5,7\mod 8$.}\endcases$$
\endpro
\par For the results related to Conjectures 4.25-4.38, see [S10].
\pro{Conjecture 4.39} Let $p$ be an odd prime.
\par $(\t{\rm i})$ If $p\e 1\mod 4$ and so $p=x^2+y^2$ with $2\nmid x$,
then
$$\aligned\sum_{k=0}^{p-1}\f{\b{-\f 14}k\b{-\f 12}k}{4^k}
&\e\Ls p3\sum_{k=0}^{p-1}\b{-\f 12}k\b{-\f 16}k2^k
 \\&\e\cases (-1)^{\f{p-1}4+\f{x+1}2}(2x-\f
p{2x})\mod{p^2}&\t{if $12\mid p-1$,}
\\2y-\f p{2y}\mod{p^2}&\t{if $12\mid p-5$.}\endcases\endaligned$$
\par $(\t{\rm ii})$ If $p\e 3\mod 4$, then
$$\sum_{k=0}^{p-1}\f{\b{-\f 14}k\b{-\f 12}k}{4^k}
\e 0\mod{p^2}\qtq{and}\sum_{k=0}^{p-1}\b{-\f 12}k\b{-\f 16}k2^k\e
0\mod p.$$
\endpro

\pro{Conjecture 4.40} Let $p$ be an odd prime. Then
$$\aligned&\sum_{k=0}^{p-1}\f{\b{-\f 14}k\b{-\f 12}k}{(-3)^k}
\e (-1)^{\f{p-1}4}\sum_{k=0}^{p-1}\f{\b{-\f 14}k\b{-\f 12}k}{81^k}
\\&\e\cases 2x-\f p{2x}\mod {p^2}&\t{if $p=x^2+y^2\e 1\mod
4$ and $2\nmid x$,}
\\0\mod p&\t{if $p\e 3\mod 4$.}\endcases\endaligned$$\endpro

\pro{Conjecture 4.41} Let $p$ be an odd prime. Then
$$\aligned&\sum_{k=0}^{p-1}\f{\b{-\f 14}k\b{-\f 12}k}{(-80)^k}
\\&\e\cases 2x-\f p{2x}\mod {p^2}&\t{if $p=x^2+y^2\e \pm 1\mod
5$ and $2\nmid x$,}
\\2y-\f p{2y}\mod {p^2}&\t{if $p=x^2+y^2\e \pm 2\mod
5$ and $2\nmid x$,}
\\0\mod p&\t{if $p\e 3\mod 4$.}\endcases\endaligned$$\endpro

\pro{Conjecture 4.42} Let $p>5$ be a prime. Then
$$\aligned&\sum_{k=0}^{p-1}\b{-\f 14}k\b{-\f 12}k2^k
\\&=\cases 2x-\f p{2x}\mod {p^2}&\t{if $p=x^2+2y^2$ with $x\e 1\mod 4$,}
\\0\mod p&\t{if $ p\e 5,7\mod 8$.}
\endcases\endaligned$$
\endpro

\pro{Conjecture 4.43} Let $p>5$ be a prime. Then
$$\aligned&\sum_{k=0}^{p-1}\b{-\f 12}k\b{-\f 13}k(-3)^k
\e\sum_{k=0}^{p-1}\f{\b{-\f 12}k\b{-\f 13}k}{(-27)^k}
\e\sum_{k=0}^{p-1}\f{\b{-\f 12}k\b{-\f 23}k}{(-4)^k}
\\&\e \Ls
p5\sum_{k=0}^{p-1}\f{\b{-\f 12}k\b{-\f 13}k}{5^k}   \e
\Ls{-1}p\sum_{k=0}^{p-1}\b{-\f 12}k\b{-\f 13}k2^k
\\&\e\cases 2A-\f p{2A}\mod{p^2}&\t{if
$p=A^2+3B^2\e 1\mod 3$ with $3\mid A-1$,}\\0\mod p&\t{if $p\e 2\mod
3$.}
\endcases\endaligned$$\endpro

\pro{Conjecture 4.44} Let $p>5$ be a prime. Then
$$\aligned&\sum_{k=0}^{p-1}\f{\b{-\f 12}k\b{-\f 13}k}{(-4)^k}
\\&=\cases \sls p5 (2A-\f p{2A})\mod{p^2}\\\q\t{if
$p=A^2+3B^2\e 1\mod 3$ with $5\mid AB$ and $3\mid A-1$,} \\ \sls p5
(A+3B-\f p{A+3B})\mod{p^2}\\\q\t{if $p=A^2+3B^2\e 1\mod 3$ with
$A/B\e -1,-2\mod 5$ and $3\mid A-1$,}
\\0\mod
p\qq\t{if $p\e 2\mod 3$.}
\endcases\endaligned$$\endpro

\pro{Conjecture 4.45} Let $p>5$ be a prime. Then
$$
\sum_{k=0}^{p-1}\b{-\f 12}k\b{-\f 16}k\Big(-\f 3{125}\Big)^k
\e\cases 2A-\f p{2A}\mod {p^2}&\t{if $p=A^2+3B^2\e 1\mod 3$,}
\\0\mod p&\t{if $p\e 2\mod 3$.}
\endcases$$
\endpro
\pro{Conjecture 4.46} Let $p>3$ be a prime. Then \par $(\t{\rm i})$
If $p\e 1\mod 4$, then
$$\sum_{k=0}^{(p-1)/2}\f{\b{2k}k^3}{64^k}\e
\Ls{3}p\sum_{k=0}^{p-1}\f{\b{2k}k\b{3k}k\b{6k}{3k}}{12^{3k}} \e
\Ls{33}p\sum_{k=0}^{p-1}\f{\b{2k}k\b{3k}k\b{6k}{3k}}{66^{3k}}
\mod{p^3}.$$
\par $(\t{\rm ii})$
If $p\e 1,2,4\mod 7$, then $$\sum_{k=0}^{(p-1)/2}\b{2k}k^3\e
\Ls{-15}p\sum_{k=0}^{p-1}\f{\b{2k}k\b{3k}k\b{6k}{3k}}{(-15)^{3k}} \e
\Ls{-255}p\sum_{k=0}^{p-1}\f{\b{2k}k\b{3k}k\b{6k}{3k}}{255^{3k}}
\mod{p^3}.$$
\par $(\t{\rm iii})$
If $p\e 1,3\mod 8$, then
$$\sum_{k=0}^{(p-1)/2}\f{\b{2k}k^3}{(-64)^k}\e
\Ls{5}p\sum_{k=0}^{p-1}\f{\b{2k}k\b{3k}k\b{6k}{3k}}{20^{3k}}\mod{p^3}.$$
\par $(\t{\rm iv})$
If $p\e 1\mod 3$, then
$$\sum_{k=0}^{(p-1)/2}\f{\b{2k}k^3}{256^k}\e
\Ls{-5}p\sum_{k=0}^{p-1}\f{\b{2k}k\b{3k}k\b{6k}{3k}}{54000^k}
\mod{p^3}.$$
\endpro
\par Let $p>3$ be a prime, and $m\in\Bbb Z_p$ with $m\not\e 0,16,64\mod
p $. From [S6] and [S7] we deduce that
$$\align \sum_{k=0}^{(p-1)/2}\f{\b{2k}k^3}{m^k}
&\e \Ls{m(m-64)}pP_{[\f p4]}\Ls{m+64}{m-64}^2
\\&\e\Ls {m(m-64)}p\Big(\sum_{x=0}^{p-1}\Ls{x^3-\f
32(3\cdot\f{m+64}{m-64}+5)x+9\cdot\f{m+64}{m-64}+7}p\Big)^2
\\&\e \Ls {m(m-16)}p\sum_{k=0}^{[p/6]}\b{2k}k\b{3k}k\b{6k}{3k}
\Ls m{(m-16)^3}^k
\\&\e \Ls {m(m-16)}p\sum_{k=0}^{p-1}\b{2k}k\b{3k}k\b{6k}{3k}
\Ls m{(m-16)^3}^k \mod p.\endalign$$

\enddocument
\bye